\documentclass[11pt]{amsart}
\usepackage{amsfonts,amsmath,amssymb}
\usepackage[all]{xy}
\usepackage{hyperref}

\begin{document}

\parindent=0pt
\parskip=6pt

\newcommand{\alg}{{\rm alg}}
\newcommand{\aff}{{\rm aff}}
\newcommand{\an}{{\rm an}}
\newcommand{\Aut}{{\rm Aut}}
\newcommand{\bc}{{\bf c}}
\newcommand{\C}{{\mathbb{C}}}
\newcommand{\Cp}{{\mathbb{C}_p}}
\newcommand{\End}{{\rm End}}
\newcommand{\F}{{\mathbb{F}}}
\newcommand{\Fo}{{\mathbb{F}^o}}
\newcommand{\BF}{{\bf F}}
\newcommand{\bF}{{\bf F}}
\newcommand{\sF}{{\sf F}}
\newcommand{\sG}{{\sf G}}
\newcommand{\G}{{\hat{\mathbb{G}}_m}}
\newcommand{\Gs}{{\mathbb{G}_{m*}}}
\newcommand{\Gt}{{\tilde{\mathbb{G}}_m}}
\newcommand{\Gl}{{\rm Gl}}
\newcommand{\Gal}{{\rm Gal}}
\newcommand{\half}{{\textstyle{\frac{1}{2}}}}
\newcommand{\Hom}{{\rm Hom}}
\newcommand{\bkappa}{{\boldsymbol{\kappa}}}
\newcommand{\hL}{{\widehat{L^\infty}}}
\newcommand{\LT}{{\rm LT}}
\newcommand{\m}{{\mathfrak{m}}}
\newcommand{\MU}{{\rm MU}}
\newcommand{\oh}{{\mathcal{O}}}
\newcommand{\ord}{{\rm ord}}
\newcommand{\Q}{{\mathbb{Q}}}
\newcommand{\iQ}{{\mathbb{Q}^\infty_p}}
\newcommand{\Qp}{{\mathbb{Q}_p}}
\newcommand{\bQp}{{\overline{\mathbb{Q}}_p}}
\newcommand{\Qs}{{\mathbb{Q}_{p*}}}
\newcommand{\QZ}{{\mathbb{Q}_p/\mathbb{Z}_p}}
\newcommand{\R}{{\mathbb{R}}}
\newcommand{\Sp}{{\rm Sp}}
\newcommand{\Spec}{{\rm Spec}}
\newcommand{\T}{{\bf T}}
\newcommand{\mT}{{\mathbb{T}}}
\newcommand{\THH}{{\rm THH}}
\newcommand{\U}{{\rm U}}
\newcommand{\bx}{{\bf x}}
\newcommand{\bxi}{{\boldsymbol{\xi}}}
\newcommand{\Z}{{\mathbb{Z}}}
\newcommand{\Zp}{{\mathbb{Z}_p}}
\newcommand{\bz}{{\boldsymbol{z}}}
\newcommand{\rig}{{\rm rigid}}
\newcommand{\ve}{{\varepsilon}}
\newcommand{\CP}{{\mathbb{C}P}}
\newcommand{\hG}{{\hat{\mathbb G}}}
\newcommand{\X}{{\mathcal X}}
\newcommand{\vkappa}{{\varkappa}}
\newcommand{\ab}{{\rm ab}}
\newcommand{\bm}{{\m_\Cp}}
\newcommand{\bQ}{{\bQp}}
\newcommand{\la}{{\langle}}
\newcommand{\ra}{{\rangle}}
\newcommand{\D}{{\mathcal{D}}}
\newcommand{\Ta}{{\bf T}}
\newcommand{\Qpa}{{\iQ}}
\newcommand{\colim}{{\rm colim}}
\newcommand{\tw}{{\widetilde{w}}}
\newcommand{\Sym}{{\sf S}}
\newcommand{\cyc}{{\rm cyc}}
\newcommand{\hoh}{{\hat{\mathfrak o}}}
\newcommand{\hoL}{{\hoh_\Li}}
\newcommand{\ohL}{{\oh_\Li}}
\newcommand{\ff}{{\mathcal F}}
\newcommand{\bexp}{{\bf exp}}
\newcommand{\blog}{{\bf log}}
\newcommand{\Spf}{{\rm Spf}}
\newcommand{\pt}{{\rm pt}}
\newcommand{\rT}{{\rm T}}
\newcommand{\tT}{{\tilde{\rm T}}}
\newcommand{\hmu}{{\hat{\mu}}}
\newcommand{\tkappa}{{\tilde{\varkappa}}}
\newcommand{\ti}{{\tilde{t}}}
\newcommand{\bep}{{\boldsymbol{\varepsilon}}}
\newcommand{\TR}{{\rm TR}}
\newcommand{\Om}{{\Omega^1}}
\newcommand{\gl}{{\rm gl}}
\newcommand{\bp}{{\rm BP}}
\newcommand{\rinf}{{\rm inf}}

\newcommand{\obQ}{{\mathcal{O}_{\overline{\mathbb{Q}}_p}}}
\newcommand{\oQ}{{\mathcal{O}_{\mathbb{Q}_p}}}
\newcommand{\tQ}{{\widetilde{\mathbb{Q}_p}}}
\newcommand{\otQ}{{\mathcal{O}_{\widetilde{\mathbb{Q}_p}}}}
\newcommand{\oL}{{\mathcal{O}_L}}
\newcommand{\Qi}{{{\mathbb{Q}^\infty_p}}}
\newcommand{\Qih}{{\mathbb{Q}^{\hat{\infty}}_p}}
\newcommand{\Li}{{L^{\infty}}}
\newcommand{\tL}{{\widetilde{L}}}
\newcommand{\otL}{{\mathcal{O}_{\widetilde{L}}}}
\newcommand{\Lih}{{L^{\hat{\infty}}}}
\newcommand{\oLih}{{\mathcal{O}_{L^{\hat{\infty}}}}}
\newcommand{\oLi}{{\mathcal{O}_{L^{\infty}}}}
\newcommand{\oQi}{{\mathcal{O}_{\mathbb{Q}^\infty_p}}}
\newcommand{\oQih}{{\mathcal{O}_{\mathbb{Q}^{\hat{\infty}}_p}}}

\newcommand{\ie}{\textit{i}.\textit{e}.}
\newcommand{\eg}{\textit{e}.\textit{g}.}

\newcommand{\bt}{{\boldsymbol{\eta}}}
\newcommand{\bbt}{{\overline{\boldsymbol{\eta}}}}
\newcommand{\tbt}{{\tilde{\boldsymbol{\eta}}}}
\newcommand{\bmu}{{\boldsymbol{\mu}}}
\newcommand{\bk}{{\boldsymbol{\kappa}}}
\newcommand{\tk}{{\overline{\boldsymbol{\kappa}}}}
\newcommand{\oC}{{\mathcal{O}_\Cp}}
\newcommand{\bg}{{\boldsymbol{\gamma}}}
\newcommand{\oP}{{\mathcal{O}_P}}
\newcommand{\hp}{{\hat{p}}}

\title{Complex orientations for $\THH$ of some perfectoid fields}

\author{Jack Morava}

\address{Department of Mathematics, The Johns Hopkins University,
Baltimore, Maryland 21218} 

\email{jack@math.jhu.edu}


\date{4 March 2020}

\begin{abstract}{This sketch argues that work of Hesselholt \cite{12} on the 
topological Hochschild homology of $\Cp$ extends, using work of Scholze 
and others \cite{28}, to define complex orientations for a version of
topological Hochschild homology for rings of integers in a natural class
of generalized cyclotomic perfectoid fields; and that the resulting spectra 
provide geometrically interesting targets for analogs of the Chern character, 
defined for certain integral lifts \cite{22} of the extraordinary $K$-functors
of chromatic homotopy theory.}\end{abstract}

\maketitle 

{\bf \S I \; Introduction and Preliminaries}

{\bf 1.1} About fifty years ago SP Novikov called attention to the relevance
of one-dimensional formal groups in algebraic topology, and since then
understanding the resulting link between homotopy theory and arithmetic 
geometry (following Quillen) has become one of the deepest topics 
in mathematics.

Recent powerful new ideas \cite{4,5,6} in $p$-adic Hodge theory,
applied to the (now classical) Lubin-Tate generalized cyclotomic closure 
$\Lih$ \cite{16} of a local number field $L$, define complex orientations (ie 
ring homomorphisms
\[
\MU^* \to \THH^*(\oLih,\Zp) \cong \Sym^*_{\oLih}(T\Omega^1_\oQ(\oLih))
\]
from the Lazard-Quillen complex cobordism ring) for the $p$-adic topological 
Hochschild homology spectra for the valuation rings of such fields. These 
orientations have interesting connections on one hand to chromatic homotopy 
theory, and on the other to the $p$-adic Fourier theory \cite{27} of Schneider 
and Teitelbaum. 
 
To present an accessible account of these connections requires considerable 
background, reviewed in this section after a sketch of the organization of 
this paper. The interests of workers in modern stable homotopy theory have 
a great deal of overlap with current work in higher local number theory, but 
the languages of these fields have diverged since the days of Cartan's
seminars. The review below summarize material useful in both areas; technical 
terms used informally in the outline immediately below will be defined more 
precisely in that review. 

{\bf 1.2 Organization} This paper is organized as follows. \S 2 interprets the 
group $\bmu$ of $p$-power roots of unity (for example in a $p$-adic field big 
enough to contain all such), as an analog in $p$-adic algebraic $K$-theory of 
the classical Hopf line bundle $\eta \to \CP^\infty$ in complex topological
$K$-theory. Building on Hesselholt's pioneering 2006 work on the Dennis trace 
\[
\tau_\Cp : k_*^\alg(\Cp,\Zp) \; \to \; \THH_*(\oC,\Zp)
\]
(for the {\bf non}-periodic $p$-adic algebraic $K$-theory of the completion 
$\Cp$ of $\bQp$), we use $k_\alg(\Cp)^*(B\bmu,\Zp)$ to construct orientation 
classes $\bg_P$ for $p$-adic $\THH$ of the valuation ring $\oP$ of a 
perfectoid subfield $P$ of $\Cp$ (containing the completion $\Qih$ of the 
field of $p$-power roots of unity over $\Qp$)\begin{footnote}{\noindent 
Unless otherwise noted, we take the prime $p$ to be odd. Some notation may 
be abbreviated, as specified below, to simplify iterated subscripts and 
similar decorations.}\end{footnote}.

The main concern of this paper is the Galois structure of Hesselholt and 
Madsen's $p$-adic $\THH$ \cite{10} of the valuation ring $\oLih$ of the 
(perfectoid) completion of a maximal {\bf totally ramified} Abelian extension 
$\Li$ of $L$, discussed in \S 3. Such questions are now accessible through 
recent extensions \cite{6} (Th 6.1) of B\"okstedt periodicity, which reduce 
this problem to that of the structure of (the $p$-adic completion of) the 
module of K\"ahler differentials of $\oLih$ over $\Zp$.

In \S 3.3, work of Fontaine from 1982 is extended to identify the Tate module 
\[
T \Omega^1_\oQ(\oLi) \cong (\pi_0 \D_L)^{-1} \oLih \otimes_\oL T_L
\]
of such objects; for example 
\[
T \Omega^1_\oQ(\oQi) \cong p_0^{-1} \oQih \otimes_\oQ T_\Qp
\]
($T_L$ being the Tate module [\S 3.1] of a Lubin-Tate group for $L$). 
Here $\D_L$ is Dedekind's different ideal of $L$ over $\Qp$, and $\pi_0$ is 
a certain torsion point of a Lubin-Tate group. A Jacobson-Zariski exact 
sequence
\[
0 \to \Omega^1_\oQ(\oQi) \otimes_\oQi \oLi \to \Omega^1_\oQ(\oLi) \to 
\Omega^1_\oQi(\oLi) = 0
\]
$p$-completes to an isomorphism 
\[
0 \to p_0^{-1} \oLih \otimes_\oQ T_\Qp \to (\pi_0 \D_L)^{-1} \oLih 
\otimes_\oL T_L \to 0
\]
defined by multiplication by an element $p_0(\pi_0 \D_L)^{-1} 
\Omega_\partial(L)\in \oLih$ closely related to the Schneider-Teitelbaum 
period of the Lubin-Tate group of $L$. The Galois behavior of the unit 
$\Omega_\partial(L)$ is key to relations between the Lubin-Tate group of 
$L$ and the complex orienation constructed for $\THH(\oLih,\Zp)$.

The final section dicusses possible applications of this construction
in chromatic homotopy theory, and an appendix summarizes some properties of 
lifts of the vector-space valued functors $K(n)$ to cohomology theories 
valued in modules over local number rings. \bigskip

\noindent
{\bf 1.3 Background}\bigskip

\noindent
{\bf 1.3.1} {\sc Conventions from local classfield theory} 

To avoid overuse of the letter $K$, in this paper $L \supset \Qp$ will denote
a locally compact topological field of degree $[L:\Qp] = n < \infty$, with a
topology defined by a discrete valuation associated to a homomorphism
\[
\ord_p : L^\times \to e^{-1}\Z \subset \Q 
\]
(normalized by $\ord_p(p) = 1$). We have 
\[
L \supset \oL \supset \m_L = (\pi_L) 
\]
with $\oL/\m_L := k_L \cong \F_q$ with $q = p^f$; this should not 
lead to confusion with $K$-theory notation. The `uniformizing element' 
$\pi_L$ generating the maximal ideal $\m_L$ of the (local) valuation ring 
$\oL$ will be fixed once and for all; it satisfies an Eisenstein equation 
\[
\pi_L^e + \cdots + u \cdot p = 0 \;
\]
with coefficients in the Witt ring $W(k_L) \subset \oL$, with $u$ a unit, 
so $\ord_p(\pi_L) =  e^{-1}$ with $n = ef$. When $L$ is unramified over $\Qp$, 
ie when $e = 1$, we will take $\pi_L = p$, and may omit the subscript on 
$\pi_L$ when $L$ is clear from context. 

We will also need {\bf non}-discretely valued topological fields, such
as the completion $\Cp$ of an algebraic closure $\bQp \supset \Qp$. The
Galois group $\Gal(\bQp/\Qp)$ acts continuously on $\Cp$ and its order
homomorphism is surjective. Similarly, $\Qi$ will denote the smallest 
subfield of $\Cp$ containg all $p$-power roots of unity, and $\Qih$ will 
be its completion; these completions are important examples of {\bf 
perfectoid} number fields \cite{27}, \ie \: $p$-adic fields $P$ complete with 
respect to a non-discrete valuation, such that the Frobenius endomorphism 
$x \mapsto x^p$ is surjective on $\oP/p\oP$. We will use $\bmu$ to denote 
the group of $p$-power roots of unity in some field of interest, fixing 
an isomorphism 
\[
\epsilon : \QZ \cong \bmu(\Qi) \;.
\]
This identifies the Tate module
\[
T\bmu := \Hom(\QZ,\bmu) 
\]
with the free rank one module $\Zp(1)$ over the $p$-adic integers, with
Galois action defined by the cyclotomic character 
\[
\chi : \Gal(\Qi/\Qp) \cong \Aut(\bmu) \cong \Zp^\times
\]
by the $p$-adic analog of the Kronecker-Weber theorem.

{\bf Example} The field $\Qp(\zeta_p) \subset \Qi$ generated over $\Qp$
by adjoining a primitive $p$th root $\zeta_p$ of unity has $f=1$ and 
$\zeta_p - 1$ as uniformizing element, satisfying 
\[
[p]_\G(x) = x^{-1}[(1 + x)^p - 1] = 0 \;,
\]
so $e = p-1$. By Lubin-Tate theory [\S 3.1.1] this is the same as the 
field $\Qp(p_0) ;= \tQ$ generated over $\Qp$ by adjoining a root $p_0$ of 
\[
[p]_{\hat{\mathbb{G}_0}}(x) = x^{p-1} + p = 0 \;:
\]
the formal groups associated to these two series are isomorphic, via 
a series $\phi(x) = x + \cdots \in \Zp[[x]]$ which sends $\zeta_p - 1$ 
to $p_0$. The Galois group of this field over $\Qp$ is cyclic, of order 
$p-1$, and the ratio $p_0^{-1}(\zeta_p - 1)$ is a unit.

{\bf 1.3.2} {\sc Algebraic $K$-theory and $\THH$} 

\noindent
{\bf 1.3.2.1} In 2003 Hesselholt and Madsen \cite{11} (\S 1.5.6) constructed, 
for a  non-Archimedean topologized field $K$ with valuation ring $A$ and 
perfect residue field, a cyclotomic spectrum $\rT(A|K)$ and a trace morphism
\[
k_\alg(K) \to \rT(A|K) \;,
\]
which refines to a map to a pro-system $\{\rm{TR}^n, \; n \geq 1\}$ of 
fixed-point sectra. These notes will unfortunately be concerned only with
the first $(n=1)$ stage of this system, which is accessible in principle by 
methods [\eg \: MacLane homology \cite{17} (\S 1)] of classical homological 
algebra, and $\rT(A|K)$ is, for such rings, a generalized Eilenberg - MacLane 
spectrum \cite{20} (\S 4). Recent work \cite{22} (\S III.5) of Scholze and 
Nikolaus has clarified, among other things, the $\mathcal{E}_\infty$ properties
of this construction, and we will be interested in the $p$-adic completion 
\[
k_\alg(K,\Zp) := (k_\alg(K))_\hp \to (\rT(A|K))_\hp =: \THH(A,\Zp)
\] 
of this trace map, defined by the functor
\[
X \mapsto X_\hp := [M(-1,\QZ),X]
\]
of maps from a suitable Moore spectrum \cite{12} (\S 2.3). 

{\bf Example} The stable $p$-adic homotopy group $\pi^S_*(B\QZ,\Zp)$, \ie  
\[
[M(\QZ,-1),[S^0,B\QZ]]_* = [M(\QZ,*-1),B\QZ] = 
\]
\[
H^1(M(\QZ,*-1),\QZ) = \Hom(\QZ,\QZ) = T(\QZ)
\] 
of $B\QZ$ vanishes unless $* = 2$, when it is the Tate module $\Zp(1)$. 
[In fact the $p$-adic completions \cite{33} (Ex. 2 p 43) of $B\QZ$ and 
$\CP^\infty$ are homotopy equivalent as spaces. 

In 2006 Hesselholt \cite{12}, greatly extending the early periodicity theorem 
of B\"okstedt \cite{7},\cite{13} (Intro), calculated the $p$-adic completion
\[
\xymatrix{
k^\alg_*(\Cp,\Zp) \ar[d]^\cong \ar[r]^{\tau_\Cp} & \THH_*(\oC,\Zp) 
\ar[d]^\cong\\
\Zp[\beta_\Cp] \ar[r] & \oh_\Cp[\gamma_p] }
\]
of this trace on the homotopy groups of these spectra, for $K = \Cp$. He
identified $\THH^*(\oC,\Zp)$ as the symmetric algebra over $\oC$ on the 
Tate module of the $p$-adic K\"ahler differentials $\Omega^1_\oQ(\oC)$, 
and showed that the class $\beta_\Cp$ (with Bockstein image $\zeta_p :=
\epsilon(p^{-1})$ in $k^\alg_1(\Cp)$) maps to $(\zeta_p - 1)$ times a 
generator
\[
\gamma_p := (\zeta_p - 1)^{-1}d \log \epsilon 
\]
of this Tate module. For reasons explained below, it will be convenient
to define a variantly normalized generator
\[
\gamma_\Cp := p_0^{-1}d \log \epsilon
\]
such that $\beta_\Cp \mapsto p_0 \cdot \gamma_\Cp$.

Hesselholt further showed that (absolute) $p$-adic Galois group acts as 
(graded) ring automorphisms of $k_\alg^*(\Cp,\Zp)$ through its abelian 
quotient
\[
\Gal(\bQ/\Qp) \to \Gal(\Qih/\Qp) \cong \Zp^\times \cong C_{p-1} \times \Zp
\]
via the cyclotomic character $\sigma(\beta_\Cp) = \chi(\sigma) \cdot 
\beta_\Cp$, compatibly with the action of $\sigma$ on $p_0 \gamma_\Cp$ by
projection to the cyclic group $\Gal(\Qp(p_0)/\Qp)$ on the first term in
the product, and by the identification of $\Zp$ with $(1 + p\Zp)^\times$
on the second. 

{\bf 1.3.2.2} In 2012 Bhatt \cite{4} showed that THH enjoys flat (ie fpqc) 
descent, and more recent work \cite{6} of Bhatt, Morrow, and Scholze shows 
that for a perfectoid $p$-adic field $P$,

{\bf Theorem} $\THH^*(\oP,\Zp)$ {\it is polynomial over $\oP$, on a single 
generator of degree two: more precisely, 
\[
\THH^*(\oP,\Zp) \cong \Sym^*_\oP(\THH_2(\oP,\Zp))
\]
is isomorphic (as $\Gal(P/\Qp)$-module) to the symmetric algebra on the Tate 
module 
\[
\THH_2(\oP,\Zp) \cong T\Omega^1_\oQ(\oP) := \Hom(\QZ,\Omega^1_\oQ(\oP))
\]
of the K\"ahler differentials of $\oP$ over $\oQ$.}

It is this result that that makes this paper possible; it will 
be applied below to generalized cyclotomic fields, whose Galois theory is
relatively well-understood. Perhaps we should explain  here that, for a 
morphism $A \to B$ of commutative rings, we write 
\[
d : B \to \Omega^1_A(B)
\]
for the universal $A$-module homomorphism satisfying $d(bc) = bdc + cdb$
(on the grounds that this notation may be easier than $\Omega^1_{B/A}$ to 
read when $A$ and $B$ are highly sub- or superscripted). 

\newpage

\begin{center}{\bf \S II The $p$-adic Hopf line bundle}\end{center}\bigskip

\noindent
The complex Hopf line bundle $\eta_\C \in k(\C P^\infty)^0(B\mT)$ 
(classified by a map $\C P^\infty \to B\mT \to B\U$) has a $p$-adic analog
\[
\eta_\Qih \in k_\alg(\Qih)^0(B\bmu,\Zp)
\]
defined by $B\bmu \to B\Gl_1(\Qih) \to B\Gl^+_\infty(\Qih)$. It is 
$\Gal(\Qi/\Qp) \cong \mathcal{O}^\times_\Qp$ invariant, and, with coproduct 
defined by multiplication
\[
B\bmu \times B\bmu \to B\bmu \;,
\] 
satisfies $\Delta(\eta_\Qih) = \eta_\Qih \otimes \eta_\Qih$; its image 
$\eta_\Cp \in k_\alg(\Cp)^0(B\bmu,\Zp)$ has similar properties.

The sequence $\{p^{-n},n \geq 1\}$, regarded via $\epsilon$ as a generator of 
$\pi^S_2(B\bmu,\Zp)$, defines a homomorphism
\[
k_\alg(\Qih)^0(B\bmu,\Zp) \to k^\alg_2(\Qih,\Zp)
\]
sending $\eta_\Qih - 1$ to a lift $\beta_\Qih$ of $\beta_\Cp$. Just as 
$\beta_\Cp$ generates $k^\alg_2(\Cp,\Zp), \; \beta_\Qih$ generates 
$k^\alg_2(\Qih,\Zp)$, so both $k_\alg(\Qi,\Zp)$ and $k_\alg(\Cp,\Zp)$ are 
complex-orientable ring spectra (which, in the case of $\Cp$, has 
been known since Suslin). We use this together with flat descent to 
construct coordinates $\bk$  (ie of cohomological degree two) for the 
related topological Hochschild groups, fitting in a diagram
\[
\xymatrix{
k_\alg(\Qih)^*(B\bmu,\Zp) \ar[d] \ar[r]^-{\tau_\Qih} & {\THH(\oQih)^*(B\bmu,
\Zp)\cong \oQih[\gamma_\Qih][[\bk_\Qih]]} \ar[d] \\
k_\alg(P)^*(B\bmu,\Zp) \ar[d] \ar[r]^-{\tau_P} & {\THH(\oP)^*(B\bmu,\Zp) \cong
\oP[\gamma_P][[\bk_P]]} \ar[d] \\
k_\alg(\Cp)^*(B\bmu,\Zp) \ar[r]^-{\tau_\Cp} & {\THH(\oC)^*(B\bmu,\Zp) \cong
\oC[\gamma_\Cp][[\bk_\Cp]]} }
\]
of graded completed Hopf algebras\begin{footnote}{We thank the referee for 
noting that (following Bourbaki [Algebra, Ch II \S 11]) we are working with 
completions of graded rings, with respect to filtrations defined using this 
grading; or, alternatively, with affine formal group schemes.}\end{footnote} 
for any perfectoid field $P$ containing $\Qih$ and contained in $\Cp$. The 
construction takes several steps.

{\bf 2.1.1}  Let $\bbt_\Qih$ be a coordinate for $k_\alg(\Qih)^2(B\bmu,\Zp)$ 
such that
\[
\beta_\Qih \bbt_\Qih = \eta_\Qih - 1 \; {\rm modulo} \; \bbt_\Qih^2 \;,
\]
with image $\bbt_\Cp \in k_\alg(\Cp)^2(B\bmu,\Zp)$; then for some $a_i
\in \Zp$ we have 
\[
\eta_\Cp - 1 = \beta_\Cp \bbt_\Cp + \sum_{i \geq 1} a_i (\beta_\Cp 
\bbt_\Cp)^{i+1} \;;
\]
so $\eta_\Cp = 1 + \beta_\Cp \bt_\Cp$ with 
\[
\bt_\Cp := \bbt_\Cp + \sum_{i \geq 1} a_i \beta_\Cp^i \bbt_\Cp^{i+1} \;.
\]
We then have a graded formal group law
\[
\Delta \bt_\Cp = \bt_\Cp \otimes 1 + 1 \otimes \bt_\Cp + \beta_\Cp \cdot
\bt_\Cp \otimes \bt_\Cp
\]
\[
:= \beta_\Cp^{-1}(\beta_\Cp \bt_\Cp \otimes 1 +_\G 1 \otimes \beta_\Cp
\bt_\Cp)
\]
on $k(\Cp)^*(B\bmu,\Zp)$ (where  
\[
X +_\G Y = X + Y + XY
\]
denotes the multiplicative formal group law, with formal power series
\[
\log_\G(X) = \log(1+X), \; \exp_\G(X) = \exp(X) - 1
\]
as logarithm and exponential). It follows that 
\[
\tilde{\eta}_\Qih = 1 + \beta_\Qih \tbt_\Qih \in k_\alg(\Qih)^2(B\bmu,\Zp) \;,
\]
with $\tbt_\Qih := \bbt_\Qih + \sum_{i \geq 1} a_i \beta_\Qih^i
\bbt_\Qih^{i+1}$, maps to $\eta_\Cp$ in $k_\alg(\Cp)^*(B\bmu,\Zp)$, and to 
$\beta_\Qih$ in $k^\alg_2(\Qih,\Zp)$. 

{\bf 2.1.2} Now by Bhatt descent 
\[
0 \to \THH_2(\oQih,\Zp)/p_0 \to (\THH_2(\oQih,\Zp)/p_0)\otimes_\oQih \oC \cong
\THH_2(\oC,\Zp)/p_0
\]
is injective, so 
\[
p_0^{-1}\tau_\Qih(\beta_\Qih) =: \gamma_\Qih \in \THH_2(\oQih,\Zp)
\]
is well-defined, and maps to $\gamma_\Cp$. 
 
{\bf 2.1.3 Lemma} {\it The formal power series $p_0^{-1} \log_\G(p_0 X) := 
\log_\Gt(X)$ and $p_0^{-1} \exp_\G (p_0 X) =: \exp_\Gt(X)$ have integral 
coefficients in $\Qp(p_0) = \tQ$, so the associated groupscheme $\Gt/\Spec \;
\otL$ is of additive type.} 

{\bf Proof} These power series are mutually inverse, so it suffices
to verify the assertion for one of them. We have 
\[
p_0^{-1} \exp_\G (p_0 X) = \sum_{n \geq 1} p_0^{n-1} \frac{X^n}{n!} \in
\Qp(p_0)[[X]] \;;
\]
but
\[
\ord_p(\frac{p_0^{n-1}}{n!}) = \frac{n-1}{p-1} - \frac{n - \alpha_p(n)}{p-1}
= \frac{\alpha_p(n) - 1}{p-1}
\]
by Legendre (where $\alpha_p(n) = \sum a_k$ when $n = \sum a_k p^k$). This is 
non-negative when $n \geq 1$. $\Box$ 

{\bf Example} $2_0 = - 2$ so $\half (e^{2x} - 1) \in \Z_{(2)}[[x]]$, 
and is congruent mod two to $\sum_{n \geq 0} x^{2^n}$. 

{\bf 2.1.4 Definition} Let
\[
\bc_\Cp := (p_0 \gamma_\Cp)^{-1} \log_\G (p_0 \gamma_\Cp \tau_\Cp(\bt_\Cp)) \in
\THH(\oC)^2(B\bmu,\Zp) \;.
\]
By the lemma, this analog of the Chern class is a power series with 
$p$-adically integral coefficients, so
\[
\tau_\Cp(\eta_\Cp) = 1 + \exp_\G(p_0\gamma_\Cp \bc_\Cp) := 1 + p_0\gamma_\Cp 
\bk_\Cp \;,
\]
with $\bk_\Cp \in \THH(\oC)^2(B\bmu,\Zp)$ a coordinate for a formal
group law with comultiplication
\[
\Delta \bk_\Cp = \bk_\Cp \otimes 1 + 1 \otimes \bk_\Cp + p_0 \gamma_\Cp
\bk_\Cp \otimes \bk_\Cp \;.
\]
A completely analogous construction defines lifts $\bc_\Qih$ and 
$\bk_\Qih$ of $\bc_\Cp$ and $\bk_\Cp$ such that
\[
\tau_\Qih(\tilde{\eta}_\Qih) = 1 + p_0\gamma_\Qih \bk_\Qih \in \THH(\oQih)^2(B\bmu,\Zp) \;, 
\]
yielding a formal group law 
\[
\Delta \bk_\Qih = \sum a_{i,j} \bk_\Qih^i \otimes \bk_\Qih^j
\]
on $\THH(\oQih)^2(B\bmu,\Zp)$. The inclusion $\Qih \to \Cp$ defines an 
an injective map on THH, so
\[
\Delta \bk_\Qih = \bk_\Qih \otimes 1 + 1 \otimes \bk_\Qih + p_0 \gamma_\Qih
\bk_\Qih \otimes \bk_\Qih \;,
\]
and thus similarly for any perfectoid field $P$ sandwiched between them. 

{\bf Corollary} {\it For such fields, the parameter $\kappa_P$ (such that 
$\tau_P(\eta_P - 1) = p_0\kappa_P$) defines a formal group law
\[
\Spf \; \THH(\oP)^0(B\bmu,\Zp) \cong \Gt \times_\otQ \oP
\]
of additive type.} 

{\bf 2.2} It follows from Hesselholt's results that $\sigma \in \Gal(P/\Qp)$
acts, as in \S 1.3.2,1, on $\gamma_P \bk_P := \kappa_P$ as multiplication by 
the cyclotomic character $\chi(\sigma)$ via the map to $\Gal(\Qi/\Qp)$. 
\bigskip

\noindent
Inverting the Bott-Thomason class $\beta_P \in k^\alg_2(P,\Zp)$
defines a 2-periodic $p$-adic algebraic $K$-theory functor 
\[
K_\alg(P)^*(-,\Zp) = \beta_P^{-1}[X,k_\alg(P,\Zp)]_{-*} 
\]
which evidently maps to $(p_0 \gamma_P)^{-1}\THH(\oP)^*(-,\Zp)$; but since
$\ord_p(p_0) > 0$, this localization defines a version 
\[
K_\alg(P)^*(-,\Zp) \to H^*(-,\gamma_P^{-1}\THH^*(\oP,\Qp))
\]
of the classical Chern-Dold character. There is a large literature
(\eg \: Segal \cite{29}, Snaith \cite{31}, Boyer, Lawson, Lima-Filho, 
Mann and Michelson \cite{8}, Totaro \cite{34,35} \dots) on related 
integrality questions.

\newpage

\begin{center}{\bf \S III Generalized cyclotomic fields}\end{center}

Any $L \supset \Qp$ as in \S 1.3.1 admits a maximal totally ramified Abelian 
extension $L^\infty \supset L$, with an Artin reciprocity homomorphism 
\cite{30}
\[
\xymatrix{
\Gal(\bQp/L) \ar[r] & \widehat{\Z} \times \Gal(L^\infty/L) \ar[r] & \Gal
(L^\infty/L) \ar[r]^{\vkappa_L} & \oL^\times 
\;,}
\]
where $\widehat{\Z} \cong \Gal(L^{nr}/L) \cong \Gal(\overline{k}_L/k_L)$,
which classifies unramified extensions, has been supressed; they will play
no part in this paper. The completion $\Lih$ of $\Li$ in $\Cp$
is perfectoid; see \cite{23} (\S 1.4.17) or \cite{38} (ex 2.0.4, ex 2.1.1) 
for its tilt, along with much more information\begin{footnote}
{To a topologist it is tempting to call these `chromatic' fields. There are
interesting analogies with the Alexander cover of a link complement
\dots}\end{footnote}.

We argue below that the graded formal groups $\THH(\oLih)^*(B\bmu,\Zp)$ 
defined above have a natural interpretation in terms of the $p$-adic Fourier 
theory of Schneider and Teitelbaum, as a rigid analytic version 
\cite{26}(\S 5 intro)
\[
\xymatrix{
\LT^\rig_L \ar[d] \ar[r]^{\ve_L} & \LT_L \ar[d] \\
\Spec \; \Cp \ar[r] & \Spec \; \oh_L }
\]
of the Lubin-Tate formal groups $\LT_L$ used to construct these 
extensions, mapped by an analog of the exponential map of classical 
Lie theory (defined only in a neighborhood of the origin). 

{\bf 3.1.1} The construction of $\LT_L$ is elegant and in some sense quite
elementary, but it depends (up to a canonical isomorphism) on some choices, 
\ie \;\; of a uniformizing element $\pi_L$ as in \S 1.3.1, as well as an element 
$[\pi_L](T) \in \oh_L[[T]]$ equal to $\pi_L T$ modulo terms of higher order, 
and congruent modulo $\m_L$ to $T^q$. In the following we will assume that
$\LT_L$ is special in the sense of Lang, ie that
\[
[\pi_L](T) = \pi_L T + T^q \;;
\]
this implies that $T \mapsto \omega T$, for $\omega \in W(k_L)^\times \subset
\oL^\times$, is an automorphism of $\LT_L$, which simplifies issues of grading
in homotopy theory. The resulting formal group  
\[
X,Y \mapsto F_L(X,Y) := X +_L Y \in \oh_L[[X,Y]]
\]
is in fact a formal $\oh_L$-module, endowed with an isomorphism
\[
a \mapsto [a]_L : \oL \to \End_\oC(\LT_L) \;.
\]
The group 
\[
\Hom_{\oL,c}(\oL[[T]],\oC) := \LT_L(\m_\Cp) = (\m_\Cp,+_L)
\]
of points of $\LT_L$ (defined by continous homomorphisms) is isomorphic 
(modulo $\Q$-vector spaces) to $(\QZ)^n$, and its Tate module 
\[
T_L := \Hom(\QZ,\LT_L(\m_\Cp))
\]
is free of rank one over $\End(\LT_L) \cong \oL$. Adjoining the torsion
points of $\LT_L(\m_\Cp)$ to $L$ defines the extension $\Li \supset L$; the
Galois group $\Gal(\bQ/L)$ acts on these torsion points, defining a 
reciprocity map
\[
\Gal(\bQ/L) \to \Gal(\Li/L) \cong \Aut_\oL(T_L) \cong \oL^\times \;.
\]
If $L_0 \subset L_1$ is Galois, then the diagram [33]
\[
\xymatrix{
1 \ar[r] & \Gal(\bQ/L_1) \ar[d]^{\vkappa_{L_1}} \ar[r] & \Gal(\bQ/L_0)
\ar[d]^{\vkappa_{L_0}} \ar[r] & \Gal(L_1/L_0) \ar[d] \ar[r] & 1 \\
{} & (L_1^\times) \:\hat{} \ar[r]^{N^1_0} & (L_0^\times) \: \hat{}
\ar[r] & \Gal(L_1/L_0)^\ab }
\]
(with $N^1_0 = N^{L_1}_{L_0}$ the norm, and profinite completion denoted by
a caret) commutes. In particular $\vkappa_\Qp(\sigma) = 
N^L_\Qp(\vkappa_L(\sigma))$ if $\sigma \in \Gal(\bQp/L)$.

Reducing modulo $\m_L$ defines a formal group law on $k_L[[T]]$ of height $n$,
and thus an embedding of $\oh_L^\times$ in the automorphism group 
$\oh_D^\times$ of the reduction as a kind of maximal torus. 

{\bf 3.1.2} A formal group law (\eg \; $\LT_L$) over a torsion-free ring 
(\eg \; $\oL$) has a unique logarithm and exponential, \eg 
\[
\log_L(T), \; \exp_L(T) \in L[[T]]
\]
such that 
\[
X +_L Y = \exp_L(\log_L(X) + \log_L(Y)) \in \oL[[X,Y]] \;.
\]
The formal multiplicative group of \S 2.1.2 is one classical example, and
Honda's logarithm \cite{14}
\[
\log_\pi(T) := \sum_{n \geq 0} \pi^{-n}T^{q^n} \in L[[T]]
\]
is another. Lemma 2.1.3 generalizes to special Lubin-Tate groups as
follows: 
 
{\bf Definition} Let $\pi_0 \in \LT_L(\m_\Cp)$ be a primitive $[\pi]$-torsion
point, ie a generator of the cyclic group of points satisfying 
\[
\pi_0^{q-1} + \pi = 0 \;,
\]
and let $\tL := L(\pi_0)$; then $\Gal(\tL/L) \cong k^\times_L \subset 
W(k_L)^\times$ (using Teichm\"uller representatives).

{\bf Proposition} {\it The formal group law
\[
X +_\tL Y := \pi_0^{-1}(\pi_0 X +_L \pi_0 Y) \in \otL[[X,Y]]
\]
is of additive type, ie with $[p]_\tL(T) \equiv 0$ modulo $\pi_0$.}

{\bf Proof} Evidently $[\pi]_\tL(T) \cong \pi(T - T^q) \cong 0$ modulo
$\pi_0$. But $[\pi]_L$ satisfies an Eisenstein equation
\[
E_L([\pi]_L) = [\pi]_L^e \; +_L \; \cdots \; +_L \; [u]_L \circ [p]_L \; = \; 0
\]
(with coefficients from $W(k_L)$, and $u$ a unit) in $\End_(\LT_L)$, so, 
similarly,
\[
[\pi]_\tL^e \; +_\tL \; \cdots \; +_\tL \; [u]_\tL \circ [p]_\tL \; = \; 0
\]
in the endomorphisms of $F_\tL$. But $[\pi]_\tL \equiv 0$ mod
$\pi_0$, so $[p]_\tL \equiv 0$ mod $\pi_0$, as well. $\Box$

{\bf Example} Honda's logarithm is $p$-typical, so its renormalization
\[
\pi_0^{-1}\log_\pi(\pi_0 T) = \sum_{n \geq 0} \pm \pi_0^{q^n-1}\pi^{-n} 
T^{q^n} 
\]
has coefficients in $\oh_L$, with
\[
\ord_p(\pi_0^{q^n-1}\pi^{-n}) = e^{-1}(\frac{q^n-1}{q-1} - n) = e^{-1}(1 +  
\cdots + q^{n-1} - n) \geq 0 \;,
\]
which goes to $\infty$ as $n$ does, making it a rigid analytic function.

Note that the completed Hopf $\otL$-algebra defined by $F_\tL$ is the 
pushforward, under the homomorphism defined by 
\[
T \mapsto \pi_0T : \oL[[T]] \to \otL[[T]] \;,
\]
of that defined by $F_L$; but (because this map does not preserve the
coordinate) it is not a morphism of formal group laws.

{\bf 3.2.1} A Lubin-Tate group has an associated $p$-divisible group
with Cartier dual $\Hom(\LT_L,\G)$; early work of Katz \cite{14} (\S 3)
identifies its Tate module
\[
T^\vee_L := \{ \beta(T) \in (1 + T \oh_\Cp[[T]])^\times \;|\: \beta
(X +_L Y) = \beta(X) \cdot \beta(Y) \}
\]
as free of rank one over $\oL$. A homotopy theorist will recognize this
as a specialization to Lubin-Tate groups of Ravenel and Wilson's almost
simultaneous description \cite{24} of the Hopf algebra $\MU_*(\CP^\infty)$
representing the Cartier dual of the universal $p$-divisible group: the
canonical inclusions $\CP^i \subset \CP^\infty$ define bordism classes $b_i$
with
\[
b(T) = 1 + \sum_{i>0} b_iT^i \;.
\]
satisfying the relation $b(X +_\MU Y) = b(X) \cdot b(Y)$ under the
Pontryagin product.

Classical Fourier analysis identifies the dual $V^* := \Hom_\R(V,\R)$ of 
a finite-dimensional vector space $V$ with its character group $V^\X := 
\Hom_c(V,\mT)$ by 
\[
V^* \ni \bxi \mapsto [\bx \mapsto \exp(i\bxi(\bx))] \in V^\chi \;.
\]
In 2001 Schneider and Teitelbaum \cite{27} (\S 2) defined a $p$-adic analog 
of the Pontryagin dual of a free $\oh_L$-module $M$ as the rigid analytic 
group $M^\X$ of locally analytic characters
\[
\theta : M \to (1 + \m_\Cp)^\times \;,
\]
and showed that the map
\[
\LT_L(\m_\Cp) \ni \alpha \mapsto [\beta \mapsto [a \mapsto
\beta([a]_L(\alpha))]] \in  \Hom_{\oh_L}(T^\vee_L,\oh_L^\X) \cong
(T^\vee_L)^\X
\]
is an isomorphism of $\Gal(\bQp/L)$-modules. The inverse of this equivalence
defines the morphism
\[
\ve_L : \LT^\rig_L := T^{\vee \X}_L \to \LT_L
\]
of group-valued functors cited above, represented by a homomorphism of
(completed) Hopf $\oh_L$-algebras from $\oh_L[[T]]$ to locally analytic
functions from $\m_L$ to $\Cp$, sending $T$ to
\[
\ve_L(T) = \exp_\G(\Omega(L)\log_L(T)) \in L_\an(\oh_L,\Cp) \;;
\]
where $\Omega(L) \in \hoL$, with $\ord_p(\Omega(L)) = (p-1)^{-1} -
e^{-1}(q -1)^{-1}$, is a `period' of the formal group $\LT_L$,
with the remarkable property that
\[
\Omega(L)^{\sigma - 1} = \{\frac{\vkappa_\Qp}{\vkappa_L}\}(\sigma) \in 
\oh_L^\times
\]
for $\sigma \in \Gal(\Li/L)$ \cite{3},\cite{27}{lemma 3.4},\cite{28} (\S 
6.2.3 Prop 6.4). This can be reformulated as the assertion that 
\[
H^1_c(\Gal(\Li/L),\oL^\times) \ni \{\frac{\vkappa_\Qp}{\vkappa_L}\}
\mapsto 0 \in H^1_c(\Gal(\Li/L),(\Lih)^\times) \;,
\]
or as the 
\
{\bf 3.2.2 Proposition} {\it The diagram (in which the vertical arrows 
multiply $T$ by the indicated element) 
\[
\xymatrix{
L[[T]] \ar[d]^{\pi_0} \ar[r]^{\ve_L} & \Lih[[T]] \ar[d]^{p_0} \\
\tL[[T]] \ar[r]^{\ve^0_L} & \Lih[[T]] }
\]
commutes; where 
\[
\ve^0_L(T) := \exp_\Gt(\Omega_0(L)\log_\tL(T)) = \Omega_0(L)T + \cdots \in 
\oLih[[T]] \;,  
\]
and $\Omega_0(L) := p_0 \pi_0^{-1} \Omega(L)$ is a unit. Moreover, the map 
$T \mapsto \ve^0_L(T)$ is equivariant, in the sense that 
\[
[\vkappa_\Qp(\sigma)]_\Gt^{-1} \circ \ve^0_L \circ [\vkappa_L(\sigma)]_\tL
= \sigma(\ve^0_L) \;,
\]
with respect to the action of $\sigma \in \Gal(L^\infty/L)$ on the 
coefficients of $\ve_L^0(T)$.} 

{\bf Proof} We have $\log_\tL([\vkappa]_L(\sigma)]_\tL(T)) = \vkappa_L
(\sigma)\log_\tL(T)$ and $\Omega(L) \vkappa_L(\sigma) = \vkappa_\Qp(\sigma)
\sigma(\Omega(L))$; while $\exp_\Gt(\vkappa_\Qp(\sigma) \cdots) = [\vkappa_\Qp
(\sigma)]_\Gt(\cdots)$. $\Box$ 

{\bf 3.3.1} In 1982 Fontaine \cite{10} (Th 1) defined, for a Lubin-Tate group 
of $L$, a homomorphism
\[
\xi_L :\bQ \otimes_\oL T_L \to \Omega^1_\oL(\oh_\bQ)
\]
of modules over the twisted group ring\begin{footnote}{with finite sums
$\sum a_\sigma \cdot \sigma, \; a_\sigma \in \oLih, \; \sigma \in 
\Gal(\Li/L)$ as elements, and multiplication $(a \cdot \sigma)(b \cdot \tau)
= a\sigma(b) \cdot \sigma \tau$}\end{footnote} $\obQ \la \Gal(\bQ/L)\ra$: if 
$h = \{h_n\} \in T_L$ with $[\pi^n]_L(h_n) = 0$, and $\bQ \ni \alpha = 
\alpha_n \pi^n$ with $\alpha_n \in \obQ$ integral, then
\[
\xi_L(\alpha \otimes h) := \alpha_n \cdot h^*_n(d\log_L(T)) \in
\Omega^1_\oL(\oh_\bQ)
\]
is well-defined; where $h_n \in \LT_L(\m_\bQ)$, and $d\log_L(T) = 
\log'_L(T) \cdot dT$. His Theorem-1 then proves the exactness of 
the sequence
\[
\xymatrix{
0 \ar[r] & \pi_0^{-1}\obQ \otimes_\oL T_L \ar[r] & \bQ \otimes_\oL T_L 
\ar[r]^{\xi_L} & \Omega^1_\oL(\oh_\bQ) \ar[r] & 0} \;;
\]
and taking Tate modules defines an isomorphism
\[
T\Omega^1_\oQ(\oh_\bQ) \cong p_0^{-1} \oh_\Cp \otimes_\oQ T_\Qp
\]
of $\Gal(\bQ/\Qp)$-modules. This argument generalizes, without significant 
change, to imply the existence of an exact sequence
\[
\xymatrix{
0 \ar[r] & \pi_0^{-1}\oLi \otimes_\oL T_L \ar[r] & \Li \otimes_\oL T_L 
\ar[r]^{\bxi_L} & \Omega^1_\oL(\oLi) \ar[r] & 0}
\]
of $\Gal(\Li/L)$-modules, and hence an isomorphism
\[
\Omega^1_\oL(\oLi) \cong (\Li/\pi_0^{-1}\oLi) \otimes_\oL T_L \;.
\] \

{\bf 3.3.2} Now a sequence $\Qp \subset L' \subset L \subset L'' \subset \Cp$ 
of extensions implies \cite{10} (\S 2.4 lemma 2) a monomorphism
\[
0 \to \Omega^1_{\oL'}(\oL) \to \Omega^1_{\oL'}(\oL'') \;,
\]
and since $\oLi$ is flat over $\oL$, the Jacobi-Zariski exact
sequence
\[
0 \to \Omega^1_\oQ(\oL) \otimes_\oL \oh_\Li \to \Omega^1_\oQ(\oLi) \to
\Omega^1_\oL(\oLi) \to 0
\]
implies that
\[
T\Omega^1_\oQ(\oLi) \cong (\pi_0 \D_L)^{-1} \oLih \otimes_\oL T_L
\]
as $\Gal(L^\infty/\Qp)$-modules; where $\Omega^1_\oQ(\oL) \cong \oL/\D_L\oL 
\cong \D_L^{-1}/\oL$, the Dedekind different being the inverse of the 
fractional ideal $\Hom_\Zp(\oL,\Zp)$ defined by the trace from $\oL$ to 
$\oQ$). 

Similarly, since $\bQ$ is flat over $\iQ$, the monomorphism
\[
0 \to \Omega^1_\oQ(\oQi) \otimes_{\oQi} \obQ \to \Omega^1_\oQ(\oh_\bQ) 
\]
implies a commutative diagram
\[
\xymatrix{
\Omega^1_\oQ(\oQi) \otimes_\oQi \oh_\bQ \ar[d]^\cong \ar[r] &
\Omega^1_\oQ(\obQ) \ar[d]^\cong \ar[r] & \Omega^1_{\oQi}(\obQ) \ar[r]
& 0 \\
\bQ/p^{-1}_0 \obQ \ar[r]^\cong & \bQ/p^{-1}_0 \obQ \;;}
\]
but then $\Omega^1_{\oQi}(\obQ) = 0$, and hence (since
\[
0 \to \Omega^1_{\oQi}(\oLi) \to \Omega^1_{\oQi}(\obQ) = 0
\]
is injective), that $\Omega^1_\oQi(\oLi) = 0$.

{\bf Corollary} {\it The exact sequence
\[
0 \to \Omega^1_\oQ(\oQi) \otimes_\oQi \oLi \to \Omega^1_\oQ(\oLi) \to
\Omega^1_\oQi(\oLi) = 0
\]
implies an isomorphism
\[
\phi:  p_0^{-1} \oLih \otimes_\oQ T_\Qp \cong (\pi_0 \D_L)^{-1} \oLih 
\otimes_\oL T_L
\]
of Tate modules over $\oLih \la \Gal(L^\infty/\Qp) \ra$.} \bigskip

\noindent
{\bf Proposition} {\it Let $t_\Qp, \; t_L$ generate the Tate modules $T_\Qp, 
\; T_L$, and let $\sigma \in \Gal(\Li/L)$; then 
\[
t^{-1}_L \phi(t_\Qp) := \Omega_\partial(L) \in \oLih
\]
also satisfies  
\[
\Omega_\partial(L)^{\sigma - 1} =  \{\frac{\vkappa_\Qp}{\vkappa_L}\}(\sigma) 
\;, \]
and hence $\D_L \Omega_\partial(L) = u \cdot \Omega(L)$ for some unit 
$u = u(L:\Qp) \in \oL^\times$.}

With $t_\Qp$ as in \S 1.3, we can choose $t_L$ so $\Omega^\partial_L \equiv
1$ mod $\m_\Lih$. The Galois action on THH of the Teichm\"uller units in
$W(k_L)$ is then consistent with the topological grading.

\begin{center}{\sc \S IV Applications and speculations}\end{center}

Writing $\ve^\partial_L$ for the variant of $\ve^0_L$ defined by replacing 
$\Omega_0(L)$ by $\Omega^\partial_0(L)$ defines a composition 
\[
\xymatrix{
\Spf \; \THH(\oLih)^0(B\bmu,\Zp) \ar[d] \ar[r]^-{\D_L} & \Gt \times_\otL \oLih
\ar[r]^-{\ve^\partial_L} & F_\tL \ar[r]^{\pi_0} & \LT_L \ar[d] \\
\Spec \; \oLih \ar[rrr] & {} & {} & \Spec \; \oL}
\]
of affine group schemes, with the top right morphism as in \S 3.1.2 and
the top left morphism defined (following \S 2.1.4) by $T \mapsto \D_L T$. 
By \S 3.2.2, the morphism across the top takes the natural action of 
$\Gal(\Li/L)$ on the left to the action of $\oh_L^\times$ by formal 
group automorphisms on the right, compatible with Artin reciprocity, and 
it seems natural to conjecture that when $L/\Qp$ is Galois, this extends
to equivariance with respect to an action of $\Gal(\Li/\Qp)$; see Appendix
III of \cite{37}, and \cite{9}.\bigskip

\noindent
The morphism $\chi_L$ from the Lazard ring to $\oL$ which classifies $\LT_L$
lifts to a graded Hirzebruch genus
\[
\MU^* \ni [M] \mapsto \chi_L[M] \cdot u^{\dim_\C M} \in \oL[u] \;,
\]
defining a commutative diagram 
\[
\xymatrix{ 
\MU^* \ar[dr]^{\tau^L_\MU} \ar[r]^{\chi_L}  & \oL[u] := k^*(L) \ar[d]^{u 
\mapsto \pi_0 \D_L \gamma_\Lih} \\
{} & \THH(\oLih,\Zp) \cong \oLih[\gamma_\Lih] }
\]
of ring homomorphisms. The appendix below summarizes a construction for weakly 
commutative complex-oriented cohomology theories with
\[
\Spf \; k(L)^0(\C P^\infty) \cong \LT_L \;.
\]
Hirzebruch's work from the 60s, interpreting multiplicative natural 
transformations 
\[
[\MU^* \to E^*] \in E^0(\MU) \cong E^0(B\U)
\]
of cohomology theories in terms of the Thom isomorphism and symmetric
functions, defines a lift of this diagram to a diagram of multiplicative
natural transformations between $\Z_2$-graded (because $\Omega^\partial_0(L)$ 
is only almost (\ie mod $\m_\Lih)$) equal to 1. \bigskip

\begin{center}{\sc Acknowledgements and thanks}\end{center}\bigskip

\noindent
This paper has roots in memorable conversations with
B\"okstedt and Waldhausen, in Bielefeld during the Chernobyl weekend: they
had just invented topological Hochschild homology, and although they could see
glimmers of connection with chromatic homotopy theory \cite{7,21}, it was
Lars Hesselholt's talk \cite{13} at the 2015 Oxford Clay symposium that really 
opened the door. I am deeply indebted to him for his generous help and patient
willingness since then to endure iterated attempts to misunderstand
his work. I hope I have not continued to garble it, and I want to thank
an extremely perceptive and insightful referee for help in clarifying 
a very rough early version of this paper.\bigskip

\noindent
I also owe thanks to Matthias Strauch for calling my attention to the relevance
of $p$-adic Fourier theory for this project, to Peter Schneider for his course
notes, to Jacob Lurie, Akhil Mathew and Chuck Weibel for conversations about 
the descent properties of THH, to Andy Baker for arithmetic and infinite 
loop-space counsel, and to Peter Scholze for sharing his insights and results, 
especially as discussed in \S 1.3.2. I am indebted as well to the 
Mittag-Leffler Institute for the opportunity for more conversations with 
Hesselholt, and to Mona Merling, Apurv Nakade, and Xiyuan Wang at JHU for 
their helpful attention. 
\bigskip

\noindent
Finally I want to thank the mathematics department at UIUC for organizing this
conference, and beyond that, to acknowledge deep lifetime debts to both Paul
Goerss and Matt Ando, two subtle, quiet enormous figures who loom, almost
invisibly, over so much of the work that has made the research world we live
in today. \bigskip

{\bf Appendix:} {\sc Integral lifts of $K(n)$ parametrized by local fields} 
\bigskip

\noindent
For any $n \in {\mathbb N}$ and prime $p$ (or, alternately, for every field
with $q = p^n$ elements), there is a functor
\[
K(n) : ({\rm Spaces}) \to (\F_q-{\rm modules})
\]
(roughly, the residue field at a prime of the sphere spectrum) with
\[
K(n)^*(B\mT) = \F_q[u^{\pm 1}][[T]]
\]
a graded complete Hopf algebra or formal group with addition $+_{K(n)}$,
with $|u| = +2, \; |T| = -2$, such that the $p$-fold multiplication
map is represented by
\[
[p]_{K(n)}(T) = u^{q-1}T^q \;.
\]
Araki's generators for $\bp_*(\pt) = \Z_p[\dots,v_i,\dots]$ satisfy
\[
[p]_\bp(\bz) \; = \; \sum_\bp v_i \bz^{p^i} = p\bz +_\bp v_1 \bz^p +_\bp
\dots
\]
($i \geq 0, \; v_0 = p$), and there is a Baas-Sullivan quotient \cite{22} :
the totalization of a suitable Koszul complex, or an iterated homotopy cofiber,
with
\[
\MU_*(\pt) \to \bp_*(\pt) \to \Z_p[u] \;,
\]
such that $v_i \mapsto 0, \; i \neq 0,n, \;; v_n \mapsto u^{q-1}$ classifying
a graded group law $+_{k(n)}$ with
\[
[p]_{k(n)}(T) = pT \; +_{k(n)} \; u^{q-1} \bz^{p^n} \;,
\]
ie $+_{k(n)} \equiv +_{K(n)} \; {\rm mod} \; p$. Hazewinkel's functional
equation  implies that
\[
T_0 +_{k(n)} T_1 = u^{-1}\exp_{k(n)}(\log_{k(n)}(uT_0) + \log_{k(n)}(uT_1)) 
\in \Zp[u][[T_0,T_1]]
\]
with
\[
\log_{k(n)}(uT) := uT + \sum_{k>0} \prod_{1 \leq i \leq k} (1-p^{q^i-1})^{-1} 
\cdot p^{-k} (uT)^{q^k} \in \Q_p[[uT]] \;.
\]
\bigskip

\noindent
Tensoring this Baas-Sullivan theory with $W(\F_q)$ defines a cohomology
theory $k(L_0)$ (where $L_0 := W(\F_q) \otimes \Q$ is the unique
{\bf un}ramified extension of degree $n$ of $\Qp$), such that
\[
k(L_0)^*(\pt) = W(\F_q)[u] \to \F_q[u] \to \F_q[u^{\pm 1}] = K(n)^*(\pt)
\]
defines a nice integral lift to a connective version of $K(n)$. More
generally, for any $L \subset \bQ$ with $n = [L:\Qp]$, there is a connective
spectrum $k(L)$ with $k(L)^*(\pt) = \oh_L[u]$, and such that $k(L)^*(B\mT)$
is a Lubin-Tate formal group for $L$; eg if $n=1$, $k(\Qp)$ is the $p$-adic
completion of connective classical topological $K$-theory (associated to the
multiplicative formal group). \bigskip

\noindent
Lubin-Tate groups of local number fields parametrize good integral lifts of
$K(n)$. In particular, the group $\oh_L^\times$ of units acts as
stable multiplicative automorphisms of $k(L)^*(\C P^\infty)$, eg with 
$\alpha \in \oh_L^\times$ sending $u$ to $\alpha \cdot u$ (generalizing 
the action of $\Zp^\times$ by $p$-adic Adams operations on classical 
($p$-completed) topological $K$-theory). It seems natural to think of 
these lifts as indexed by maximal toruses in the unit group $D_n^\times$ 
of a division algebra with center $\Qp$ and Brauer-Hasse invariant $1/n$. 
Recent work \cite{19} of Hopkins and Lurie, using the modern theory of Thom 
spectra, has changed the geography of this subject: in particular, it raises 
the question of possible Azumaya multiplications on such lifts, which might
support lifts of the $\Gal(\Li/L)$-action discussed above, to an action
compatible with some such multiplicative structure. 

\bibliographystyle{amsplain}

\end{document}